\documentclass[11pt]{article}
\usepackage{a4wide}
\usepackage{amsmath}

\usepackage{amsfonts}
\usepackage{amssymb}
\usepackage{euscript}
\newenvironment{proof}{\noindent {\it Proof.~~}\ }{\  \rule{1mm}{2mm}\medskip}
\newenvironment{example}{\noindent {\bf Example.~~}\ }{\  \rule{1mm}{2mm}\medskip}

\newenvironment{proof*}{\noindent {\it Proof.~~}\ }{}

\newtheorem{theorem}{Theorem}
\newtheorem{lemma}[theorem]{Lemma}
\newtheorem{corollary}[theorem]{Corollary}
\newtheorem{proposition}[theorem]{Proposition}

\newenvironment{claim}{\noindent {\bf Claim.~~}\it }{\medskip}
\def\Z{\mathbb Z}
\def\zp{{\mathbb Z}/p{\mathbb Z}}
\def\zn{{\mathbb Z}/n{\mathbb Z}}

\newcommand{\subgp}[1]{\langle{#1}\rangle}
\newcommand{\ot}{\overline{T}}
\newcommand{\ost}{\overline{S+T}}
\begin{document}
\title{On the critical pair theory in abelian groups : Beyond Chowla's Theorem}

\author{ Yahya O. Hamidoune\thanks{Universit\'e Pierre et Marie Curie, Paris {\tt yha@ccr.jussieu.fr} }
\and Oriol Serra\thanks{Universitat Polit\`ecnica de Catalunya,
Barcelona {\tt oserra@mat.upc.es}. Supported by the Spanish Research
Council under project
  MTM2005-08990-C02-01 and by the Catalan Research Council under
  project 2005SGR00256} \and Gilles Z\'emor\thanks{
Institut de Math\'ematiques de Bordeaux, Universit\'e Bordeaux 1,
{\tt zemor@math.u-bordeaux1.fr}
}}

\date{}
\maketitle

\begin{abstract}
We obtain critical pair theorems for subsets $S$ and $T$ of an
abelian group such that $|S+T|\leq |S|+|T|$. We generalize some
results of Chowla, Vosper, Kemperman and a more recent result due
to R{\o}dseth and one of the authors.
\end{abstract}

\section{Introduction}
Let  $S$ and $T$ be nonempty  subsets of  $\zp$. The Cauchy-Davenport Theorem
\cite{CAU,DAV} states that
$$|S+T|\geq \min (p,|S|+|T|-1).$$

The Cauchy-Davenport Theorem was
generalized to abelian groups by several authors including Mann
\cite{MANA} and Kneser \cite{KN3}. The first  generalization to cyclic
groups is due to Chowla \cite{CHO}: it states
\begin{theorem}\label{th:chowla}
Let $S,T $ be nonempty subsets of $\zn$ such $0\in S$.
Assume that every element of $S\setminus \{0\}$ has order exactly $n$.
Then $|T+S|\geq \min(n, |S|+|T|-1).$
\end{theorem}

Subsets achieving equality in an additive theorem are known as critical
pairs of the theorem. One may easily
check that the only interesting critical pairs for the Cauchy-Davenport Theorem arise
when $|S|,|T|\ge 2$ and $|S+T|\leq p-2$. Under these assumptions
 Vosper's Theorem \cite{VOS} states that $|S+T|=|S|+|T|$ unless both $S$
and $T$ are   arithmetic
progressions with a common difference. This statement determines the
critical pairs of the Cauchy-Davenport Theorem.

Generalizing Vosper's Theorem to arbitrary abelian groups requires a
lot of care. The importance of this question was mentioned by Kneser
in \cite{KN3}. Motivated by Kneser's work,  Kemperman proposed in
\cite{KEM} a recursive procedure which generalizes Vosper's Theorem
to abelian groups. The main tools used by Kemperman are basic
transformations introduced by Cauchy, Davenport and Dyson
\cite{Nat}. One of the results obtained by Kemperman is the
following:

 \begin{theorem}[Kemperman, \cite{KEM}]\label{th:KK}
 Let $G$ be an finite abelian group and let $S,T$ be  subsets of $G$ such that
$|S|\geq 2$, $|T|\geq 2$ and $|S+T|= |S|+|T|-1\le p-2,$
where $p$ is the smallest prime divisor
 of $|G|$. Then $S$ and $T$ are  arithmetic progressions with the same
difference.
 \end{theorem}

Note that the existence of a short direct proof for this result is
unlikely since the statement contains Vosper's Theorem. This result
has been recently extended to non abelian groups by K\'arolyi
\cite{karolyi} and independently by one of the authors \cite[Theorem
3.2]{HALG}.

By using the  additive transformations mentioned above, R{\o}dseth
and one of the authors recently characterized the critical pairs of
Vosper's Theorem \cite{HR} :

\begin{theorem}\label{th:HR}
 Let $S, T$ be subsets of a group of prime order $\zp$,
with $|T|\ge 3$ and $|S|\ge 4$ such that
  $$|S+T|= |S|+|T|\le p-4.$$
Then $S$ and $T$ are included in arithmetic progressions with the same
difference and of respective lengths $|S|+1$ and $|T|+1$.
\end{theorem}

There are several methods  currently available in additive theory.
One of them is based on Fourier analysis.
Examples of applications of this method can be found
the monographs of Freiman \cite{FRE} and Tao and Vu \cite{TV}, or in the
papers by Deshouillers and Freiman \cite{DF}, and by Green and Ruzsa
\cite{GR}. Another powerful tool is  the polynomial method introduced
 by Alon, Nathanson and  Ruzsa \cite{ANR}.
K\'arolyi recently \cite{KAR} used this method to obtain a
remarkable critical pair theorem for restricted sums.

In this paper we obtain improvements of some of the above results
using  the {\it isoperimetric} method.
This  method has been used to generalize addition theorems to non abelian
 groups in some papers including \cite{ZEM,HAST,HALG,AOY}.
 It also derives additive inequalities, mainly from the
 structure of the {\em $k$-atoms} of a set. If $S$
 is a generating subset containing $0$ of an abelian group~$G$,
 a set $A$ is called a $k$-atom of $S$ if
 it is of minimum cardinality among subsets $X$ such that
 $|X|\geq k$, $|X+S|\leq |G|-k$, and
 $|X+S|-|X|$ is of minimum possible cardinality
 (see Section~\ref{isop} for detailed definitions).
It is proved in \cite{HCAY} that any
$1$-atom containing $0$ is a subgroup. This result implies easily
Mann's generalization of the Cauchy-Davenport Theorem.
The structure of $2$-atoms has proved more difficult to describe
but potentially gives stronger results: $2$-atoms have been used
in \cite{HA00,HP} to derive critical pair results.  In groups of
prime order, the description of  $2$-atoms was completed by two of the
present authors in \cite{SZ}. Atoms of higher order were used in
\cite{GOY}
to classify  sets $S,T\subset \zp$ with $|S+T|\leq |S|+|T|+1$.

In the present paper we first study the structure of $2$-atoms in general
abelian groups. Our main result in the first part of this paper
is Theorem \ref{th:main}: broadly speaking
it states that, under some technical conditions that will be shown to
be quite tight, $2$-atoms  have cardinality~$2$ or are subgroups.
In the rest of the paper we apply this fact to obtain critical
pair results.

We shall first obtain a critical pair result for Chowla's
Theorem \ref{th:chowla}
which reduces to Vosper's Theorem if $n$ is a prime. To be precise, we
will actually be dealing with a strengthened version of Theorem
\ref{th:chowla} (Corollary \ref{chowla}) that only requires the order
of every element of $S\setminus \{0\}$ to exceed $|S|-1$ rather than
to equal $n$. We call this requirement a {\it weak Chowla condition}.
The description of the corresponding sets $S$ and $T$ are obtained
in Theorem \ref{th:kempermannis} and Corollary~\ref{th:kemperman+0}.

We then move on to give a description of subsets
$S,T$, with $|S+T|\leq |S|+|T|,$ in arbitrary abelian groups
provided  $S$ contains no element of order less than $|S|+1$
(another weak Chowla condition). We
show that,  if the abelian group has no subgroups of order $2$ or $3$,
then $S$ and $T$ are made up of arithmetic progressions with at
most one missing element and periodic subsets with at most one missing
element, see Theorems \ref{th:kempermannis+1} and
 \ref{th:kemperman+1}. This last result is a generalization
to abelian groups of Theorem \ref{th:HR} of R{\o}dseth and one of the
authors, since it reduces to it when the group is of prime order.

The paper is organized as follows:
Section \ref{isop} gives some preliminary results and
Section~\ref{sec:chowla} uses them to derive a solution to the critical pair
problem for Chowla's Theorem and its strengthened version.
Section \ref{sec:fainting} works out some tools necessary to
 Section \ref{sec:2atoms} which is devoted to the description of
$2$-atoms.
Sections \ref{sec:small} and \ref{sec:quasi} make up more preliminary
material for section~\ref{sec:chowlavosper} which derives the
generalization to abelian groups of Theorem \ref{th:HR}.

\section{Isoperimetric tools}\label{isop}


In this section we recall known results on isoperimetric numbers
of subsets in finite abelian groups and derive some consequences relevant
to us later on. Our prime objects of concern are the $2$-atoms of
a subset: we shall see that they are either subgroups or Sidon
sets and, in the last case, they have the largest possible
isoperimetric numbers.

Let  $S$ be a
subset of a finite abelian group such that $0\in S$. Denote by $\subgp{S}$ the
subgroup generated by $S$. For a positive integer $k$, we
shall say that $S$ is $k$-{\em separable} if there exists
$X\subset \subgp{S}$ such that $|X|\geq k$ and  $|X+S|\leq |\subgp{S}|-k$.

Suppose that $S$ is $k$-{separable}. The
$k$-th {\em isoperimetric number} of $S$  is then defined by
\begin{equation}  \label{eq:kappa}
\kappa _k (S)=\min  \{|X+S|-|X|\  \Big| \ \
X\subset\subgp{S}, \ |X|\geq k \ {\rm and }\ |X+S|\leq |\subgp{S}|-k\}.
\end{equation}

 For a $k$-separable
set $S$, a subset $X$ achieving the above minimum is called a
$k$-{\em fragment} of $S$. A $k$-fragment with minimal cardinality
is called a $k$-{\em atom}.

The following easy facts will be used regularly throughout the paper:

\begin{itemize}
\item
if $S$ is $k$--separable then  $1\leq\kappa_{k-1}(S)\le
\kappa_{k}(S)$.
\item
The translate $A+g$ of a $k$--atom $A$ is also a $k$--atom.
\end{itemize}

{\bf Remark}.
Let  $0\in S$  be a $k$-separable subset of a finite abelian group such that
$|S|\geq k$.
Then  $\kappa_k (S)\leq k|S|-2k+1$.

\begin{proof}
Assume the contrary. Let $G=\subgp{S}$. Then we must have clearly
$|G|\geq 2k+\kappa _k(S)\geq k|S|+2.$ Hence $k|S|-k+1\leq|G|-k-1$.
Let $X$ be a $k$-subset of $S$ such that $0\in X$.

We have $|S+X|\leq |S|+\sum _{x\in X\setminus 0} |(S+x)\setminus S|
\leq |S|+(k-1)(|S|-1)\leq k|S|-k+1< |G|-k.$

Therefore, by (\ref{eq:kappa}), we have
$\kappa _k(S)\leq |S+X|-|X|\leq k|S|-2k+1,$ a contradiction.
\end{proof}

If $S$ is not $k$-separable, we shall put by convention
$\kappa _k(S)= k|S|-2k+1$ so as to have, for all $|S|\geq k$,
\begin{equation}  \label{eq:maxkappa}
\kappa _k(S)\leq k|S|-2k+1.
\end{equation}

The definition of a $k$-atom implies the following lemma:
\begin{lemma}\label{lem:2surjk} Let $0\in S$ be a $k$-separable subset
of a finite abelian group. Let $A$ be a $k$-atom and suppose that
$|A|>k$.
Then, for each $a\in A$ and $s\in S$ we have
$$
(A\setminus \{ a\})+S=A+S=A+(S\setminus \{ s\}).
$$
\end{lemma}

\begin{proof} Let $A'=A\setminus \{ a\}$ and suppose that $|A'+S|<|A+S|$.
Then $|A'+S|-|A'|\le |A+S|-1-|A'|=|A+S|-|A|$ contradicting the
minimality of $A$. In other words, no element $x$ in $S+A$ can be
uniquely written as $x=s+a$, $s\in S$ and $a\in A$. This means that
$A+S=A+(S\setminus \{ s\})$ for each $s\in S$.
\end{proof}

Next we recall:

\begin{lemma}[\cite{HALG}]\label{lem:dual}
Let  $0\in S$  be a $k$-separable  subset of a finite
abelian group $G$. Let  $F$  be a  $k$-fragment of $S$ and $g\in
\langle S\rangle$.
Then $g-F$ and   $\langle S\rangle\setminus (F+S)$   are $k$-fragments of $-S.$
Moreover $\kappa _k(-S)=\kappa _k(S).$
\end{lemma}

The following is a particularly useful property of $k$-atoms.

\begin{lemma}[The intersection property
  \cite{HALG}]\label{lem:intersection}
Let  $0\in S$  be a  $k$-separable subset of a finite abelian group $G$.
Let $A$  be a  $k$-atom of $S$. Let  $F$  be a $k$-fragment of $S$
such that  $ A \not\subset F $.  Then   $|A \cap F|\leq k-1 .$
\end{lemma}

The intersection property implies easily the following description
of $1$-atoms.

 \begin{corollary}[\cite{HCAY}] \label{Cay}
Let  $0\in S\neq \langle S\rangle$
be a  subset of a finite abelian group $G$.
Let $A$  be a $1$-atom of $S$ such that $0\in A$.
 Then  $A$ is the subgroup  generated by $S\cap A$. In particular $\kappa _1(S)$ is a multiple of $|A|$.
\end{corollary}

From these early results we can derive the following generalization
of Chowla's Theorem:

\begin{corollary}\label{chowla}  Let  $0\in S$  be a generating subset of a
finite abelian group $G$  such that  the order of every element of
$S\setminus \{0\}$ is at least  $|S|-1$. Then $\kappa_1(S)=|S|-1$.

In particular, for every nonempty subset $X\subset G$, we have
$$
|X+S|\ge \min\{ |G|, |X|+|S|-1\}.
$$
\end{corollary}

\begin{proof}  If $S$ is not $1$-separable, then by definition we have
$S=\subgp{S}$ and by the convention preceding \eqref{eq:maxkappa} we
have $\kappa_1(S)= |S|-1$. Suppose therefore that $S$ is
$1$-separable. Let $A$ be a $1$-atom of $S$ containing $0$. By
Corollary \ref{Cay}, $A$ is the subgroup of $G$  generated by $S\cap
A$ and $\kappa_1(S)$ is a multiple of $|A|$. If $S\cap A = \{0\}$
then it follows that $A$  is the null subgroup and we have $\kappa
_1(S)=|S+A|-|A|=|S|-1$. If $S\cap A\neq \{0\}$ then by the
hypothesis on the order of the elements of $S$ we have
$\kappa_1(S)\geq |A|\geq |S|-1$ which implies $\kappa_1(S)=|S|-1$ by
\eqref{eq:maxkappa}.

The last inequality in the statement is a direct consequence of
the definition of $\kappa_1$.
\end{proof}

Recall that a subset $X$ of an abelian group is a Sidon set if no
two pairs of (not necessarily distinct) elements in $X$ have the
same sum. In particular $|S\cap (S+x)|\le 1$ for each $x$.

\begin{corollary}\label{cor:2atom}
 Let  $0\in S$  be a $k$-separable subset of a finite abelian
 group $G$. Let $A$  be a $k$-atom of $S$ such that $0\in A$, and suppose that $p\geq k$ where $p$ is the smallest prime divisor of $|G|$.
 Then either $A$ is a subgroup of $G$ or $|A\cap (x+A)|\leq k-1$ for every $x\in G$, $x\neq 0$.
In particular a $2$-atom of a $2$-separable set is either a
subgroup  or a Sidon set.
\end{corollary}

 \begin{proof}
Without loss of generality we may suppose $\subgp{S}=G$.
The double inequality $k\leq |A\cap (x+A)|<|A|$ is forbidden by   Lemma \ref{lem:intersection} because $x+A$ is also a $k$-atom of  $S$. Suppose that  there is $x\in   G$, $x\neq 0$, such that $A=A+x$. Then we have $A =  A+\subgp{x}$: hence $A\cap (a+A) \supset a+\subgp{x}$ for every   $a\in A$. Since $|\subgp{x}|\geq p\geq k$, Lemma   \ref{lem:intersection} implies that we have $A=a+A$ for every $a\in  A$ and $A$ is a subgroup.
\end{proof}

\begin{lemma}\label{Sidon0}
Let $0\in S$ be a generating set of a finite abelian group $G$
of cardinality $|S|\ge 3$.
Assume that $|(S+g)\cap S|\leq 2$ for all $g\in G\setminus \{
0\}$. Then $\kappa_1(S)=|S|-1$.
 In particular $\kappa_1(X)=|X|-1$, if $X$ is a Sidon set containing $0$.
\end{lemma}

\begin{proof} Suppose on the contrary that $\kappa_1(S)\le |S|-2$.
Then $S$ is $1$-separable.
Let $0\in A$ be a $1$-atom of $S$.
Then $A$ is a nonnull subgroup of $G$ and $\kappa_1 (S)$ is a
multiple of $|A|$.

Let $a\in A\setminus \{ 0\}$. We have $|S|+|A|-2\ge |S+A|\ge
|S\cup (S+a)|\ge 2|S|-2$ which implies
$|S|\le |A|\le \kappa_1(S)$, contradicting \eqref{eq:maxkappa}.
\end{proof}

The next  result  determines the second isoperimetric number of
Sidon sets. In what follows we use the following notation. Given a
subgroup $H$ of $G$, by the decomposition of a subset $S\subset G$
modulo $H$ we mean the minimal partition of $S$ into nonempty
subsets, each one contained in a single coset of $H$.

\begin{lemma}\label{lem:sidon} Let $0\in S$ be a subset of a
  finite abelian group
 with  $|S|\ge 3$. If $S$ is a Sidon set then $\kappa_2(S)=2|S|-3$.
\end{lemma}

\begin{proof}
Let $G=\subgp{S}$. Suppose $S$ is $2$--separable, otherwise the result
follows by the convention preceding \eqref{eq:maxkappa}.

Suppose against the lemma that $\kappa_2(S)\le 2|S|-4$. Let $A$ be a
$2$-atom of $S$ containing $0$.  We have  $|A|\ge 3,$ since
otherwise $2+\kappa _2(S)=|S+A|\ge |S|+(|S|-1),$ a contradiction. By
Corollary \ref{cor:2atom}, $A$ is either  a Sidon set or a subgroup.
We have
$$
|A|+2|S|-4\geq |A|+\kappa_2(S)=|S+A| = |\cup_{a\in A}(S+a)|\ge
|S|+(|S|-1)+(|S|-2)= 3|S|-3,
$$
which gives $|A|\geq |S|+1\ge 4.$

If $A$ is a Sidon set, then
$$
|S+A| = |\cup_{s\in S}(s+A)|\ge |A|+(|A|-1)+(|A|-2)\ge |A|+2|S|-1,
$$
a contradiction.

Suppose that $A$ is a subgroup. Then $\kappa_2(S)$ is a multiple of
$|A|$. In particular, $|A|\le 2|S|-4$. But then, since $|A|\ge 4$,
$$
|S+A| = |\cup_{a\in A}(S+a)|\ge |S|+(|S|-1)+(|S|-2)+(|S|-3)\ge
|A|+2|S|-2.
$$
again a contradiction.
\end{proof}

The following corollary is a result obtained in a more general
context  in \cite{HA00}. The simple proof given here is
 similar to a proof given in \cite{SZ}.

\begin{corollary}\label{cor:m+3}
Let $S$ be a generating set of the finite abelian group $G$ with
$0\in S$, $|S|\ge 3$ and $\kappa_2 (S)=|S|+m$, $m\geq -1$. Let $0\in
A$ be a $2$-atom of $S$ which is not a subgroup of $G$. Then $|A|\le
m+3$.
\end{corollary}

\begin{proof}
Suppose on the contrary that $|A|\ge m+4$. By Corollary
\ref{cor:2atom}, $A$ is a Sidon set of $G$. By Lemmas \ref{Sidon0}
and \ref{lem:sidon}, we have $\kappa_1 (A) = |A|-1$ and
$\kappa_2(A)=2|A|-3$.

If $A$ generates $G$ then $2|A|-3=\kappa_2 (A)\le |S+A|-|S|=|A|+m$,
a contradiction. Therefore we may assume that $A$ generates a proper
subgroup $Q$ of $G$. Let  $S=S_1\cup \cdots \cup S_j$, where $j\ge
2$, be the decomposition of $S$   modulo $Q$. We may assume  that
$|S_1+A|\leq  \cdots \leq |S_j+A|$ and, by translating $S$, that
$0\in S_1$.

If $|S_2+A|\le |Q|-1$ then,
$$
2|A|-2=2\kappa_1 (A)\le |A+S_1|-|S_1| +|A+S_2|-|S_2|\leq
|A+S|-|S|=|A|+m,
$$
against our assumption. Therefore we may assume that $S'+A=S'+Q$
where $S'=S\setminus S_1$.

If $|S_1|=1$ then, for each $2$-subset $X$ of $Q$, we have
$|S+X|-|X|=|(S+X)\setminus (S_1+X)|\le
|S'+Q|=|S'+A|=|S+A|-|S_1+A|=|S+A|-|A|$ contradicting that $A$ is a
$2$-atom of $S$ with $|A|\geq 3$. Hence $|S_1|\ge 2$. Now if
$|S_1+A|\le |Q|-2$ then $2|A|-3=\kappa_2 (A)\le |S_1+A|-|S_1|\le
|S+A|-|S|=|A|+m$. Hence we may assume $|S_1+A|\geq |Q|-1$.

Since $|S+Q|-1\le |S'+Q|+|A+S_1|=|A+S|\leq |G|-2,$ we must have
$S+Q\neq G.$ But in this case,
$$
|S|+m=\kappa _2(S)\leq |S+Q|-|Q|\leq |S+A|+1-|Q|= |S|+|A|+m-|Q|+1.
$$
It follows that
 $|A|\ge |Q|-1,$ which is impossible since $A$ is a Sidon set.
\end{proof}

Finally, the following lemma will be useful to us in ruling out the
possibility that a $2$-atom is a subgroup.

\begin{lemma}\label{lem:order}
Let $0\in S$ be a $2$-separable subset of a finite abelian group $G$.
Suppose $A$ is a $2$-atom of $S$ which is a subgroup of cardinality
at least $3$. Then there exists $s\in S$, $s\neq 0$, such that
the order of $s$ is not more than $\kappa_2(S)$.
\end{lemma}

\begin{proof}
 Note that if $A$ is a subgroup then $\kappa_2(S)$ is a multiple of
 $|A|$. By Lemma \ref{lem:2surjk} we have $A+S=(A\setminus\{0\})+S$
 which implies $(A\setminus\{0\})\cap (-S)\neq \emptyset$. Therefore
 there is a non-zero element $s$ of $S$ in $A$, and its order is
 not more than $|A|\leq\kappa_2(S)$.
\end{proof}

\section{Critical pairs under the weak Chowla condition.}\label{sec:chowla}

With the previous results we can already prove a critical pair
theorem  improving  on the theorems of Chowla and Vosper. We first
state its  isoperimetric version. Recall that a subset $S$ of an
abelian group $G$ is {\it periodic} if there is a nonnull subgroup
$H$ of $G$ such that $S+H=S$. In other words, $S$ is a union of
cosets of $H$.

\begin{theorem}\label{th:kempermannis}  Let $0\in S$ be a generating
$2$-separable subset of a finite abelian group $G$ such that
$\kappa _2(S)\leq |S|-1.$
Also assume  that every element of $S\setminus \{ 0\}$ has order
at least $|S|$. Then either
 $S$ is an arithmetic progression or  $S\setminus \{ 0\}$ is periodic.
\end{theorem}

\begin{proof}
By Corollary \ref{chowla} we have $\kappa_1 (S)=|S|-1$. Let $0\in
A$ be a $2$-atom of $S$. Assume $|S|\geq 3$ otherwise there is nothing
to prove. By Lemma \ref{lem:order}, the condition on the order of
elements of $S$ implies that $A$ is not a subgroup. But then
Corollary \ref{cor:m+3} implies that we have
$|A|=2$, say $A=\{0,r\}$. Assume first that $r$ generates $G$.
This forces   $S$ to be an arithmetic progression with difference
$r$. Assume now that $r$ generates a proper subgroup $H$. Let
$S=S_1\cup \ldots S_j$, $j\ge 2$ be the  decomposition of $S$
modulo $H$. We have $|S|+1= \sum_{i=1}^ j|S_i+\{0,r\}|\ge
\sum_{i=1}^ j\min\{ |H|, |S_i|+1\}$, which implies $|S_i|=|H|$ for
all but one subscript.  In particular $S\cap H=\{ 0\}$ since
otherwise $S$ contains a nonzero element with order at most
$|H|\leq |S|-1$.
\end{proof}

The above theorem will translate into a Chowla-type
characterization  of sets $S$ and $T$ with small sumset, this will be
Corollary \ref{th:kemperman+0}. The next
result is a generalization of Theorem~\ref{th:KK}.

By the {\em stabilizer} of a subset $X$ of an abelian group $G$, we mean
the set of group elements $x\in G$ such that $X+x=X$.

\begin{proposition}\label{fragments}
Let $0\in S$ be a generating subset of a finite abelian group $G$ and
let $0\in T$ be a subset of $G$. Let $Q$ denote the stabilizer of
$S\setminus \{0\}$. Suppose that
  $$|T+S|\leq |T|+|S|-1< |G|-|Q|.$$
Also assume that every element of $S^*=S\setminus \{0\}$ has order $\geq |S|$.
Let $\sigma:G\rightarrow G/Q$ denote the canonical projection. One of
the following holds:\vspace{-2mm}
\begin{itemize}
\item[(i)]  either $T\subset Q$,\vspace{-2mm}
\item[(ii)] or $\sigma (S)$ and $\sigma (T)$ are  arithmetic progressions with
      the same difference. Moreover, at most one member of the decomposition
      of $T$ modulo $Q$ is not a complete coset modulo~$Q$.
\end{itemize}
\end{proposition}

\begin{proof} Either $T=\{ 0\}$ and thus $T\subset Q$ or the conditions on $S$ imply that
$S$ is $2$--separable and  $\kappa_2 (S)\leq |S|-1$. Assume first
$Q=\{0\}$. By Theorem~\ref{th:kempermannis}, $S$ is an arithmetic
progression. It follows easily that $T$ is an arithmetic progression
with the same difference. Assume now $Q\neq \{0\}$.

We have  $|\sigma(T)+\sigma(S)|\leq |\sigma(T)|+|\sigma(S)|-1$.
Otherwise there are $|\sigma(S)|$ cosets in $\sigma (T)+\sigma(S)$ not
present in $\sigma(T)$. But all these cosets are saturated in $T+S$
(notice that $S^*$ is $Q$-periodic). It follows that
$|T+S|\ge |T|+|\sigma (S)||Q|=|T|+|S|+|Q|-1$, a contradiction.

 Moreover, the order of every element $x\in \sigma (S)\setminus \{0\}$
 is at least $\lceil |S|/|Q|\rceil =|\sigma (S)|$. Since the stabilizer
 of $\sigma(S)^*=\sigma(S^*)$ must be $\{0\}$, either $\sigma (T)=0$
 and $T\subset Q$ or
Theorem \ref{th:kempermannis} in $G/Q$ implies that $\sigma(S)$ is
an arithmetic progression. It follows now that $\sigma (T)$ is an
arithmetic progression with the same difference. Since $\sigma(T)$
contains at most a single element that is not expressible in $G/Q$
in two different ways as a sum of one element of $\sigma(S)$ and one
element of $\sigma(T)$, we deduce that at most one coset modulo $Q$
that intersects $T$ is not included in $T$.
\end{proof}

\begin{corollary}\label{th:kemperman+0}  Let $0\in S$ and $T$ be
non-empty subsets of a finite abelian group $G$. Suppose that
$$
|S+T|\leq |S|+|T|-1< |H+T|-|Q|,
$$
where  $Q$ denotes the stabilizer of $S\setminus \{0\}$ and $H$ is the
subgroup of $G$ generated by $S.$ Also assume
that every element of $S\setminus \{ 0\}$ has order at least
$|S|$. Let $T_1\cup T_2\cup \cdots \cup T_j$  be a decomposition of $T$
modulo $H$ such that
$$
|T_1+S|\leq |T_2+S|\leq \cdots \leq |T_j+S|.
$$
Then $|T_i|=|H|$ for all $i\geq 2$.
Moreover one of the following conditions holds:
\begin{itemize}
\item[(i)]
$T_1-T_1\subset Q.$
\item[(ii)] $\sigma(S)$ and $\sigma(T_1)$ are arithmetic progressions with the same difference,
where $\sigma:G\rightarrow G/Q$ denotes the canonical projection.
\end{itemize}
\end{corollary}

\begin{proof}
By Corollary \ref{chowla} we have $\kappa_1(S)=|S|-1$. If $j\geq 2$
we have $|T_2+S|=|H|$ since otherwise,
$$
2|S|-2=2\kappa_1 (S)\leq |S+T_1|-|T_1| +|S+T_2|-|T_2|\leq
|S+T|-|T|\leq |S|-1,
$$
a contradiction. Assume first that $0\in T_1$. Since $|S+T|<
|H+T|-|Q|$, we have $|S+T_1|<|H|-|Q|$. We have clearly
$$
|S+T_1|\le |S|+|T_1|-1< |H|-|Q|.
$$
 By Proposition \ref{fragments}, either $T_1\subset Q$ or
$\sigma(S)$ and $\sigma(T_1)$ are arithmetic progressions with the
same difference. Now, if $0\not\in T_1$ then the same argument gives
$T_1-T_1\subset Q$.
\end{proof}

At the heart of the proof of Theorem \ref{th:kempermannis} was the
claim that, under the right conditions, a $2$-atom containing the
zero element is of cardinality $2$ or is a subgroup. In Section
\ref{sec:2atoms} we shall find more general conditions under which
we can make the same claim. Before that we need some more tools.

\section{The fainting technique}\label{sec:fainting}

In this section we use a method developed in \cite{SZ}. The idea is
to consider the sequence of subsets $(S+A)\setminus S$,
$(S+2A)\setminus (S+A), \cdots , (S+iA)\setminus (S+(i-1)A),\cdots$
and to claim that if $A$ is a $2$-atom of $S$ of cardinality
$|A|>2$, then this sequence must decrease and faint, implying that
$S$ is a ``large'' subset of $G$.

Let  $X$ and $Y$ be subsets of an abelian group $G$.
For each integer $i\ge 0$ we
denote by $$N_i(X,Y)=(X+iY)\setminus (X+(i-1)Y),\; i>0,\;\;
N_0(X,Y)=X,$$ where $iY=\underbrace{Y+\cdots +Y}_{i}$.

In what follows we use the notation  $Y^*=Y\setminus \{ 0\}$. We start
with the two following lemmas.

\begin{lemma}  \label{commutativity}
Let $G$ be an abelian group and
let $X,Y\subset G$ with $0\in X\cap Y$.
If $N_r(X,Y)-Y^*\subset N_{r-1}(X,Y)$ for some
$r\ge 1$, then $N_i(X,Y)-Y^*\subset N_{i-1}(X,Y),$ for all $i\geq r$.
\end{lemma}

\begin{proof}
Suppose that the statement holds for all $i$, $r\le i\le j$, for some
$j\ge r$, and let $x\in N_{j+1}(X,Y)$ (if $N_{j+1}(X,Y)=\emptyset$ there
is nothing to prove.) By the definition of $N_{j+1}(X,Y)$, there
is $z\in Y^*$ such that $x-z\in N_{j}(X,Y)$. Now, for every $y\in
Y^*$, $x-y-z=(x-z)-y\in N_{j-1}(X,Y)$ which implies $x-y\in
N_{j}(X,Y)$. The result follows by induction.
\end{proof}

\begin{lemma}\label{lem:surj}
Let $0\in S$ be a $2$-separable subset of a finite abelian group $G$ and
let $0\in A$ be a $2$-atom of $S$ with cardinality $|A|\geq 3$
which is not a subgroup of $G$.
Then, denoting $A^*=A\setminus\{0\}$,
  $$S+A=S+(A\setminus \{ a\}) \mbox{ for each } a\in A
\mbox{   and } N_2(S,A)-A^*\subset N_1(S,A).$$
\end{lemma}

\begin{proof}
Without loss of generality $S$ generates $G$.
The first part of the result is just Lemma~\ref{lem:2surjk}. Now,
since $A$ is not a subgroup, we have $S+\subgp{A}\neq S+A$,
otherwise we would have $|S+\subgp{A}|-|\subgp{A}|< |S+A|-|A|$ in
contradiction with $A$  being a $2$-atom. Therefore there exists
$x\in N_2(S,A)=(S+2A)\setminus(S+A)$.
Recall that, by Lemma \ref{lem:dual}, the subset $x-A$ is a $2$-atom
of$-S$ and $G\setminus(S+A)$ is a $2$-fragment.
Observe that $x\in (x-A)\cap (G\setminus (S+A))$ and
 that $x\in N_2(S,A)$ means $x-A$ is not contained in $G\setminus
 (S+A)$:
the intersection property of $2$-atoms (Lemma \ref{lem:intersection})
implies therefore that $(x-A)\cap (G\setminus (S+A)=\{x\})$,
but this means $x-A^*\subset N_1(S,A)$.
\end{proof}

The following Lemma is a key tool for the proof of the main
result of the next section.
It says that, under some conditions, a set $X$ verifying
the statement of Lemma \ref{lem:surj} with some other set   must
be a large subset of the ground group.

\begin{lemma}[The Fainting Lemma] \label{fainting}
Let $G$ be a finite abelian group and
let $X,Y\subset G$ with
$0\in X\cap Y$ and set $m=|X+Y|-|X|-|Y|$. Assume that $Y$
generates $G$ and that
\begin{itemize}
\item[(i)]
$3\leq |Y|\le m+3$  and $\kappa_1(Y^*-y)=|Y^*|-1\ge 1$ for some $y\in Y^*$.
\item[(ii)]
$X+Y=X+(Y\setminus \{ z\})$ for each $z\in Y$   and
$N_2(X,Y)-Y^*\subset N_1(X,Y)$.
\end{itemize}
Then
$$|X|\ge |G|-\binom{m+4}{2}.$$
\end{lemma}

\begin{proof}
Since $X+Y=X+(Y\setminus \{ y\})$ for any $y\in Y$, we have $X+Y=X+Y^*$ and $X+(Y-y)=X+(Y^*-y)$.
By induction on $i$ it is seen that $X+i(Y-y)=X+(Y^*-y)+(i-1)(Y-y)=X+i(Y^*-y)$ for each $i\ge 1$. Since $0\in Y$ generates $G$, we have $G=X+n(Y-y)=X+n(Y^*-y)$ where $n=|G|$. Let
 $H$ be the subgroup of $G$ generated by $Y^*-y$. One can verify easily that  $H=n(Y^*-y)$, and hence
\begin{equation}\label{eq:f1}
X+H=G.
\end{equation}

By Lemma \ref{commutativity}, $N_2(X,Y)-Y^*\subset N_1(X,Y)$
implies
\begin{equation}\label{eq:f2}
N_{i+1}(X,Y)-Y^*\subset N_{i}(X,Y) \mbox{ for all } i\ge 1.
\end{equation}
Fix $y\in Y^*$ satisfying (i). Suppose that there is $i\ge 1$ such
that $N_{i+1}(X,Y)\neq\emptyset$ and

\begin{equation} \label{eq:cannot}
|N_{i+1}(X,Y)-(Y^*-y)|<|N_{i+1}(X,Y)|+|Y^*|-1.
\end{equation}

Since $\kappa_1 (Y^*-y)=|Y^*|-1$, the inequality (\ref{eq:cannot})
means that $N_{i+1}(X,Y)-(Y^*-y)$ is a union of cosets of the subgroup
$H$ generated by $Y^*-y$. In particular  $H\neq G$ and, by
(\ref{eq:f2}),  $N_i(X,Y) \supset(N_{i+1}(X,Y)-(Y^*-y))+y$
contains a full coset of this subgroup. However, we have
$N_i(X,Y)\cap X=\emptyset$ and, by (\ref{eq:f1}), $X+H=G$, a
contradiction.
Let $\ell$ be the largest integer for which
$N_{\ell}(X,Y)\neq\emptyset$. We have just shown that, for each $i$,
$1\le i<\ell$,
$$
|N_{i+1}(X,Y)|+|Y^*|-1\le
|N_{i+1}(X,Y)-(Y^*-y)|=|N_{i+1}(X,Y)-Y^*|\le |N_i(X,Y)|.
$$
Therefore,
\begin{equation}\label{eq:bm+4}
|G|=
|X|+\sum_{i=1}^{\ell}|N_i(X,Y)|\le|X|+\sum_{i=1}^{\ell}(|N_1(X,Y)|-(i-1)(|Y|-2)).
\end{equation}
Since $|N_1(X,Y)|=|Y|+m$ we have $|N_2(X,Y)|\le m+2$. Hence, since
$3\le |Y|\le m+3$,  the largest possible value in the right hand
side of inequality (\ref{eq:bm+4}) is taken if $|Y|=3$ and $\ell
=m+3$ giving
$$
|G|\le |X|+\binom{m+4}{2},
$$
as claimed.
\end{proof}

We finish this set of preliminary results with the following
Lemma.

\begin{lemma}\label{BigS}
Let $A$ and $S$ be subsets of a finite abelian
  group $Q$. Assume that $|A|=3$ and that for each $a\in A$ we have
$S+A=S+(A\setminus \{ a\})$. Then $3|S|\geq 2|S+A|.$
\end{lemma}

\begin{proof} Write $A=\{x,y,z\}$.  We have  $S+A=(x+S)\cup (y+S)$. It
  follows that $|S+A|=2|S|-|(x+S)\cap (y+S)|$. Furthermore we must
  have
$((x+S)\cup (y+S)) \setminus ((x+S)\cap (y+S)) \subset z+S$,
therefore $|(x+S)\cap (y+S)|\geq |S|/2$. The result now follows.
\end{proof}

\section{Description of $2$-atoms}\label{sec:2atoms}

The next theorem gives the structure of the $2$-atoms for not too
large subsets of an abelian group.

\begin{theorem}\label{th:main}
Let $G$ be a finite abelian group and let $0\in S$ be a generating $2$-separable
subset of $G$ such that $|S|\ge 3$ and $$\kappa_2(S)-|S|=m\le 4.$$
Let $A$ be a $2$-atom of $S$ containing $0$. If $|S|<
|G|-\binom{m+4}{2}$ then either $|A|=2$ or $A$ is a subgroup of
$G$.
\end{theorem}

\begin{proof}
Suppose that the conclusion of the theorem does not hold, so that
$0\in A$ is a $2$-atom of $S$ with $|A|\ge 3$ which is not a
subgroup. Then it follows from Corollary \ref{cor:2atom}  that $S$ is a Sidon
set and then, by Lemma \ref{lem:sidon}, $\kappa_2 (S)= 2|S|-3\ge |S|$.  In
particular $m\ge 0$.

By  Corollary \ref{cor:m+3} we have $|A|\leq
m+3$. Moreover,  $A^*-a$ is also a Sidon set, and Lemma \ref{Sidon0}
implies that $A$ satisfies condition $(i)$ of the Fainting Lemma. By
Lemma \ref{lem:surj}, $S$ and $A$ satisfy condition $(ii)$ of the
Fainting Lemma: therefore if $A$ generates $G$ its conclusion must
hold. In that case we have $|S|\ge |G|-\binom{m+4}{2}$ against the
hypothesis of the Theorem. Therefore $A$ must generate a proper
subgroup $Q$ of $G$. Let $S=S_1\cup S_2\cup \cdots \cup S_{t}$ be
the decomposition of $S$ modulo $Q$ and $I=\{1,\ldots, t\}$. Put
\begin{eqnarray*}
U&=&\{i\in I : |A+S_i|=|Q|\},\\V&=&\{i\in I : |A+S_i|=|Q|-1\},\\
W&=&\{i\in I : |A+S_i|\leq |Q|-2\}, \mbox{ and }\\
&&u=|U|, v=|V|, w=|W|.
\end{eqnarray*}

Since $|S+A|=|S|+|A|+m$, the decomposition  of $S+A$ modulo $Q$
gives
\begin{equation}  \label{eq:|S|+|A|+m}
|S|+|A|+m=
\sum_{i=1}^{t}|S_i+A|=|S|+\sum_{i=1}^{t}(|S_i+A|-|S_i|)\geq
|S|+\sum_{i\in V\cup W}(|S_i+A|-|S_i|)
\end{equation}

Now, as mentioned above,
$A$ is a Sidon set and by Lemma
\ref{Sidon0}, we have $\kappa_1(A)=|A|-1$. Therefore
$|S_i+A|-|S_i|\geq |A|-1$ for $i\in V$. Notice furthermore that
Lemma \ref{lem:2surjk} implies $S_i+A=S_i+A^*$,
so that $|S_i|\geq 2$ for each $i\in I$.
Therefore, since by Lemma \ref{lem:sidon}  we have $\kappa_2
(A)=2|A|-3$, we have $|S_i+A|-|S_i|\geq 2|A|-3$ for $i\in W$.
Inequality \eqref{eq:|S|+|A|+m} gives us

\begin{equation}\label{eq:m2}
4\ge m\ge v(|A|-1) + w(2|A|-3) -|A|.
\end{equation}
In particular we have
$$
w\le 2\mbox{ and } v\le 3.
$$

Now for any $i\in  I$ let us write
$$
\delta(i)=|Q|-|S_i+A|,\hspace{5mm}\mbox{and}\hspace{5mm}
|S_i+A|=|S_i|+|A|+m_i,
$$
and, for   $J\subset I$, put $\delta (J)=\sum _{i\in J}\delta(i)$.
Notice that $\delta(I)=|S+Q|-|S+A|$, $\delta(U)=0$, $\delta(V)=v$
and that we have shown that  $m_i\geq |A|-3\geq 0$ for $i\in W$.
We consider two cases.

\noindent{\it Case 1.\/} $S+Q\neq G$.

It follows that, by the minimality in the definition of a
$2$-atom, $|S+Q|-|Q| \geq |S+A|-|A|$. Therefore,

\begin{equation} \label{eq:|Q|-|A|}  \delta(V)+\delta(W)=|S+Q|-|S+A| \geq |Q|-|A|.\end{equation}

Since $A$ is a Sidon set we have $|Q|-|A|\geq 4$ (for example use
$|Q|\geq |A|(|A|+1)/2$ and rule out the case $|A|=3$ and $|Q|=6$ by
exhaustive search.) By  \eqref{eq:m2} and since  $\delta(V)=v\leq 3$
we have  $w\ge 1$ which in turn implies  $v\leq 2$. Now, by the
definition of $m_i$ and $\delta(i)$, we have
$|Q|=|S_i|+|A|+m_i+\delta(i)$ for all $i$ and inequality
\eqref{eq:|Q|-|A|} can be rewritten as

\begin{equation}\label{eq:delta(i)}
  \delta(W)\geq |Q|-|A|-v\geq |S_i|+m_i-v+\delta(i).
\end{equation}

Let $i\in W$. If $|A|\geq 4$ then  $m_i\geq |A|-3\geq 1$ so that
$|S_i|+m_i\ge 3$ and, if $|A|=3$, then Lemma~\ref{BigS}  and $m_i\ge
0$ give $3|S_i|\geq 2(|S_i+A|)\ge 2(|S_i|+|A|)$ meaning $|S_i|\geq
6$. In both cases inequality \eqref{eq:delta(i)}   gives $\delta
(W)>\delta (i)$, which implies $w\ge 2$.

Returning to (\ref{eq:m2}) it follows that $v=0$, $w=2$, $|A|=3$
and \eqref{eq:|S|+|A|+m} gives $\sum_{i\in W}m_i\leq 1.$
We may assume $W=\{1,2\}$. Since
$|S_i|\geq 4$ for $i\in W$, inequality  ~\eqref{eq:delta(i)} gives
$$
\delta(1)+\delta(2)\geq 4+\delta(i), \; i=1,2.
$$
Hence $\delta(i)\geq 4$ for $i=1,2$. Now since $m_1+m_2\leq 1$ we
have, for example, $m_1=0$. But the Fainting Lemma applied  to
$S_1$ and $A$ gives $6=\binom{m_1+4}{2}\geq |Q|-|S_1|$
contradicting $|Q|-|S_1|= |A|+\delta(1)\geq 7$.

\noindent{\it Case 2.\/} $S+Q=G$.

Now $|G|=|S+A|+\delta (V)+\delta (W)$ and the hypothesis of the
theorem reads
  \begin{equation}   \label{eq:hypothesis}
\delta(V)+\delta(W)+|A|+m\geq 1 +\binom{m+4}{2}.
 \end{equation}

Since $|A|\leq m+3$ (Corollary \ref{cor:m+3}), inequality
\eqref{eq:hypothesis}   implies $\delta(V)+\delta(W)\geq 4$,
giving $w\geq 1$ since $\delta(V)\leq 3$.

If $w=1$, say $W=\{ 1\}$, then \eqref{eq:hypothesis} translates to

\begin{equation}  \label{eq:hypothesis2}
 |Q|-|S_1|=|A|+m_1+\delta(1)\ge 1+\binom{m+4}{2}+m_1-m-\delta(V).
\end{equation}

If $m=m_1$ then we must have $\delta(V)=0$ and the right hand side
of \eqref{eq:hypothesis2} equals $1 +\binom{m_1+4}{2}$. If $m_1<m$
then $\delta(V)=v\leq 3$ and $m_1\geq 0$ imply that the right hand
side of\eqref{eq:hypothesis2} is again $\geq 1 +\binom{m_1+4}{2}$.
In both cases this contradicts the Fainting Lemma applied to $S_1$
and $A$.

If $w=2$, say $W=\{ 1,2\}$, then (\ref{eq:m2}) implies $|A|=3$,
$v=0$, $0\le m_1+m_2\le 1$ and $m\geq 3$. Then, the Fainting Lemma
applied to $S_i$ and $A$, $i=1,2$,  gives
$|A|+m_i+\delta(i)\leq\binom{m_i+4}{2}$. By adding up  the two
inequalities we get
$$
2|A|+m_1+m_2+\delta(1)+\delta(2)\leq
\binom{m_1+4}{2}+\binom{m_2+4}{2}.
$$
But the right hand side is at most $16$ while, by
\eqref{eq:hypothesis}, the left hand side is at least
$\delta(W)+|A|+m-1\geq \binom{m+4}{2}\ge 21$ a contradiction. This
completes the proof.
\end{proof}

The following example shows that the result of Theorem
\ref{th:main} does not hold anymore if $m=5$.

\begin{example}
Take $G=\Z/7\Z\times\Z/q\Z$ where $q>7$ is a prime. Consider the
sets

\begin{itemize}
\item $S=\{0,1,2,4\}\times X$ where $|X|=4$ and $X$ is a  Sidon set in $\Z/q\Z$ and
\item $A=\{0,1,3\}\times \{ 0\}$.
\end{itemize}

Then $|S+A|=|S|+|A|+5$. The group $G$ has only two proper subgroups
$H_1=\Z/7\Z\times \{ 0\}$ and $H_2=\{ 0\}\times \Z/q\Z$ and
$$|S+H_1|=|S|+|H_1|+5, |S+H_2|=4q>|S|+q+5.$$
On the other hand, if $B=\{ 0,x\}$ we have
$$|S+B|\ge \left
\{\begin{array}{ll}   |S|+8 >|S|+|B|+5, & \hbox{ if } x\in H_1 \\
|S|+12> |S|+|B|+5,
& \hbox{ if } x\not\in H_1\\
\end{array}\right.$$

The last inequality being because, for any $y\neq 0$, $|X\cup
(X+y)|\geq 7$ in $\Z/q\Z$ since $X$ is a Sidon set. Therefore
subgroups and subsets of size $2$ are not $2$--atoms of $S$.
Furthermore we have $\kappa_2(S)\geq |S|+5$ since otherwise
Theorem~\ref{th:main} would apply: therefore $A$ is a $2$-atom of
$S$  and $\kappa_2(S)= |S|+5$.
\end{example}

Finally, note that Theorem \ref{th:main} together with Lemma \ref{lem:order}
give a sufficient condition to rule out the possibility of a $2$-atom being a
subgroup.

\begin{corollary}\label{cholaatom}
Let $G$ be a finite abelian group and let $0\in S$ be a generating $2$-separable
subset of $G$ such that $|S|\ge 3$ and $$-1\leq
\kappa_2(S)-|S|=m\le 4.$$ Also assume  that every non zero element
of $S$ has order at least $|S|+m+1$.

If $|S|< |G|-\binom{m+4}{2}$ then  the $2$-atoms of $S$ have
cardinality $2$.
\end{corollary}


\section{Atoms of small sets}\label{sec:small}

We next show some results about $k$-atoms of small sets.


\begin{lemma}\label{4at}
Let $S$  be a $4$-separable generating subset
 of a finite abelian group such
that $0\in S$ and $\kappa _4(S)=|S|=3.$ Let $0\in A$ be a
$4$-atom of $S$. Then   $|A|=4$.
\end{lemma}

\begin{proof}
Let $G=\langle S\rangle$.
Suppose that $|A|>4$. We shall apply the Fainting Lemma to $A$ and $S$.

We have $\kappa_1(S^*-z)=1$, for all $z\in S^*$, since  $\kappa_1(S^*-z)>0$.

Take $z\in S^*$. By Lemma \ref{lem:2surjk} we have
$A+S=A+\{0,z\}=A\cup (A+z).$ Therefore $((A+S)\setminus A)\subset
A+z$. Therefore $N_1(A,S)-z\subset A=N_0(A,S)$. It follows that
$$N_1(A,S)-S^*\subset
A=N_0(A,S).$$ By  Lemma \ref{commutativity} we have
$N_2(A,S)-S^*\subset N_{1}(A,S)$.

Now we may apply the Fainting Lemma and obtain $|A|\ge |G|-6$. But
then $|A+S|\ge |G|-3$ contradicting that $A$ is a $4$-fragment of $S$.
\end{proof}

\begin{lemma}\label{lem:|A|<=3} Let $0\in S$ be a $3$-separable generating
subset of a finite abelian group $G$  such that $\kappa_3(S)=|S|=4$.
Assume  $\gcd (|G|,6)=1.$  Let $A$ be a $3$-atom of $S$ such that
$0\in A$. Then $|A|=3$.
\end{lemma}

\begin{proof} Suppose on the contrary that $|A|\ge 4$.
Then $A$ is not a subgroup since otherwise $|S+A|$ and $\kappa_3(A)$
are multiples of $|A|$, so that
$|A|\le \kappa_3(S)\le 4$ contradicts $\gcd (|G|, 6)=1$.
By Corollary \ref{cor:2atom} we have
\begin{equation}\label{eq:2int1}
|A\cap (A+g)|\leq 2 \mbox{ for each } g\neq 0.
\end{equation}
In particular,
$$
|S|+|A|=|S+A|=|\cup_{s\in S} (s+A)|\ge |A|+(|A|-2)+(|A|-4),
$$
which implies $|A|\le 5$.

Suppose that $A$ generates a proper subgroup $H$ of $G$ and let
$S=S_1\cup \cdots \cup S_j$ be the decomposition of $S$ modulo
$H$. By Lemma \ref{Sidon0},
$$
|S_i+A|\geq \min (|H|,|S_i|+|A|-1), \; i=1,\ldots ,j.
$$
Choose  $h\in H\setminus \{ 0\}$. Since $|S_i|\le 3$,
$$
|H|\ge |A\cup (A+h)|\ge 2|A|-2\geq |S_i|+|A|-1.
$$
Therefore, $|S+A|=\sum_{i=1}^j|S_i+A|\ge |S|+2|A|-2>|S|+|A|$, a
contradiction. Hence $\subgp{A}=G$.

In particular,  $|S+A|=|S|+|A|\le |G|-3$ implies $\kappa_2 (A)\le
\kappa_3 (A)\le |A|$. Let $0\in B$ be a $2$-atom of $A$. Since
$|A|\le |G|-7$, Theorem \ref{th:main} implies that  $B$ is a
subgroup or $|B|=2$.

Suppose that $B$ is a subgroup. Then $|B|\ge 5$. Let $A=A_1\cup
\ldots \cup A_j$ be the decomposition of $A$ modulo $B$. We have
$(j-1)|B|=|A+B|-|B|\le |A|\le 5$ which implies $j=2$, $|A|=|B|=5$
and $\kappa_2(A)=|A|$. But then $S$ is a set with smaller
cardinality than $B$ with $|S+A|-|S|=|A|$, contradicting the
minimality of the $2$-atom.

Therefore $|B|=2$. Then, using \eqref{eq:2int1}, $|A|+|B|\ge
|A+B|\ge |A|+(|A|-2)=2|A|-2$ which implies $|A|=4$.

Since $S$ is $3$-separable, we have $|G\setminus (A+S)|\geq 3$.
But we must have $|G\setminus (A+S)|>3,$ since otherwise, by Lemma
\ref{lem:dual}, $-S$ has a $3$-atom  $T$ with size $|T|=3$. This
would imply that $-T$ is a $3$-atom of $S$, a contradiction. Since
$\gcd (|G|,6)=1$, we must have $|G\setminus (A+S)|=|G|-8\geq 5.$

\begin{claim}
$S^*$ is an arithmetic progression.
\end{claim}

Let us write  $N_i=N_i(A,S)$, $i\ge 0$. Note that
$|N_1|=|S+A|-|A|=\kappa_3(S)=|S|=4$. For each subset $X\subset S$
and for each $i\geq 1$, let us denote by $N_i^X$ the set of elements
$u\in N_i$ such that $u-X\subset N_{i-1}$ and $X$ is a maximal
subset of $S$ with this property. By the definition,
$N_i^X=\emptyset$ whenever $0\in X$. Moreover, for two different
subsets $X, Y$, we have $N_i^X\cap N_i^Y=\emptyset$.

Let $X\subset S^*$, $i\geq 2$, and $u\in N_i^X$, so that $v=u-x\in
N_{i-1}$ for each $x\in X$. Let $Y$ be the subset of $S^*$ such that
$v\in N_{i-1}^{Y}$, implying $v-y\in N_{i-2}$ for each $y\in Y$.
Then $u-y=(v-y)+x\in N_{i-1}$ implies that $y\in X$. We have just
shown that:

\begin{equation}\label{eq:nx1}
\mbox{for $i\geq 2$,}\hspace{1cm}N_i^X-X\subset \cup_{Y\subset X}N_{i-1}^ Y.
\end{equation}

 By Lemma \ref{lem:2surjk}, for each $x\in S^*$, we have $A+S=A+(S\setminus \{ x\})$,
 which implies $N_1^{\{x\}}=\emptyset$.
 By \eqref{eq:nx1}, we have $N_i^{\{ x\}}=\emptyset$ for each $i\ge 1$ as well.

On the other hand, for each $x\in S^*$, inequality
\eqref{eq:2int1} implies
$$
2\le |(A+x)\setminus A|\le |N_1\setminus N_1^{S^*\setminus
\{x\}}|=4-|N_1^{S^*\setminus \{x\}}|,
$$
so that $|N_1^{S^*\setminus \{x\}}|\le 2$.

Let us now estimate $|N_2^X|$ and $|N_3^ X|$ for $X\subset S^*$.
Note that by Corollary \ref{chowla} as $\kappa_1 (Z)=|Z|-1$ for each
subset $0\in Z\subset G$ with $|Z|\leq 3$, since the order of any
nonzero element in $G$ is at least $5$. Therefore, using
\eqref{eq:nx1}, we have
$$
\mbox{ for each $2$-subset $X$ of $S^*$, }\; |N_2^X|+1\le
|N_2^X-X|\le |N_1^X|\le 2 \mbox{ and } N_3^ X=\emptyset.
$$
Since there are at most two $2$-subsets of $S^*$ for which
$|N_1^X|=2$, we have
  $$\sum_{X\subset S^*, |X|=2}|N_2^X|\le~2.$$

Since $|N_1|=4$, then $N_2^{S^*}-S^*$ cannot be a coset. Therefore,
since  $\kappa_1(S^*-s)=|S^*|-1$, we have
\begin{equation}\label{eq:ns*1}
|N_2^{S^*}|+2\le |N_2^{S^*}-S^*|\le |N_{1}|.
\end{equation}
This implies that $|N_2^{S^*}|\le 2$.

Suppose that $|N_2^{S^*}|\le 1$. Then $|N_2|=\sum_{X\subset S^*}
|N_2^ X|\le 3$ and, by applying \eqref{eq:nx1} with $i=3$ and $4$,
we get $|N_3|=|N_3^{S^*}|\le 1$ and $|N_4|=0$. Therefore
$|N_2|+|N_3|\le 4<|G|-|S+A|$. This means that $Y=A\cup N_1\cup
N_2\cup N_3\neq G$ and $Y+S=Y$, which contradicts that $S$
generates~$G$.

Suppose now that $|N_2^{S^*}|=2$. Then $|N_2^{S^*}-S^*|=|S^*|+1$
which implies that $S^*$ is an arithmetic progression. This proves
the claim. 

Now we have $S=\{ 0, a, a+d, a+2d\}$ for some $d\in G$. By
repeating the argument of the claim 
to $S-a-d$ we get that $\{-a-d,-d,d\}$ is an arithmetic progression
as well. We cannot have $-a-d=0$. Then  either $-a-d-d=2d$ or
$-a=-2d$. Hence either  $S=\{ -4d,-3d,-2d,0\}$ or  $S=\{ 0,2d,3d,
4d\}$. But then $\{ 0,d,2d\}$ is a $3$-atom of $S$. This
contradiction concludes the proof.
\end{proof}


\section{Quasi-progressions}\label{sec:quasi}

A subset $S$ of an abelian group $G$ will be called a
{\it quasi-progression of difference $r$ } if $S$ is not a
progression with difference $r$ and  if $S$ can be obtained by
deleting an element of an arithmetic progression of difference
$r$.

\begin{lemma}\label{lem:transfer}  Let $0\in S$ be a quasi-progression
with difference $r$  in the  cyclic group  $\zn$.
Suppose that $S$ generates $\zn$ and $|S|\ge 3$. Let $T\subset \zn
$ be such that $|T|\geq 3$ and
$$
|S+T| \leq |S| + |T|\leq n-4.
$$
Then one of the following conditions holds:

(i)  $T$ is either a quasi-progression with  difference  $r$
or  a progression with  difference  $r$.

(ii) $n=12$ and $T$ is a coset of order $4$.
\end{lemma}

\begin{proof}
Put $S=\{a,a+d, \cdots ,a+(j-1)d, a+(j+1)d,   \cdots ,|S|d\}.$
Observe that $S\subset \subgp{d}+a$. Since $0\in S$, we have $a\in \subgp{d}$.
Then $\zn=\subgp{S}=\subgp{d}$.
Hence without loss of generality we may assume that $j\geq \lceil| S|
/2\rceil$ and  $a=0$. Since $d$ is invertible  we
can assume it to equal $1$. Then we have
\begin{equation}\label{eq:T}
S = \{0,1,\ldots ,j-1,j+1 \cdots , |S|\}.
\end{equation}

For a subset $X\subset \zn$ let us call {\it connected components}
of $X$ the maximal arithmetic progressions with difference $1$
contained in $X$.

{\it Case 1.\/} There is a connected component $C_1$ of
$\ot=\zn\setminus T$ such
that $|C_1|\geq |S|$.

Then we clearly have $|C_1\cap (S+T)| \geq |S| -1$. Furthermore,
since $\{0,1\}\subset S$, we have $|C\cap (S+T)|\geq 1$ for every
connected component $C$ of $\ot$. Therefore $\ot$, and hence $T$,
has at most two connected components. If $\ot = C_1$  we are done.
Suppose   that $\ot = C_1 \cup C_2$ where $C_2$ is the other
component of $\ot$. Since $|T|\ge 3$, one of the two components of
$T$, say $T_1$, has cardinality $|T_1|\ge 2$. Since $S$ is a
quasi-progression, $S+T_1$ is an arithmetic progression of length
$|S|+|T_1|$. If $|C_2| \geq 2$ then we must have either
$|C_2\cap(T_1+S)|\ge 2$ or $|C_1\cap (T_1+S)|\ge |S|$, a
contradiction. Therefore we must have $|C_2| =1$, which proves the
result.

{\it Case 2.\/} For every connected component $C$ of $\ot$, $|C| <
|S|$.

It follows that every connected component $C$ of $\ost$ has size $1$.
Then  $q=|\ost|$ is  the number of connected components of $\ost$. Then
$|S+T| \geq |T| + q(j-1) \geq |T| + q(\lceil |S| /2\rceil-1)$.
Since $q\geq 4,$ we must have $q=4$, $|S|=4$ and $j=2$,
i.e. $S=\{0,1,3,4\}$.

If $U=\{u,u+1, \cdots , v\}$ is a connected component of $T$, then $v+1\in (S+T)\setminus T$.
Since $|S+T|\le |S|+4$, it follows that $T$ has exactly 4 components, and for
each such component $U$, $v+2\notin S+T$ but $v+3\in T$. Thus
$v-1\notin U$, and
each component of $T$ has exactly one element. This shows that $|T|=4$
and hence $n=12$.

 $S=\{0,1,3,4\}=\{0,1\}+\{0,3\}.$ Now $\ost$ consists of $4$
 single-element components. It follows that  $|-\ost+\{0,1\}|=|\ost|+4.$
Therefore $-\ost +S=-\ost+\{0,1\}+\{0,3\}=-\ost+\{0,1\}$.
It follows that $\ost -S$ is a union of cosets modulo the subgroup $H$
generated by $3$. Therefore $T=G\setminus (\ost -S)$ is an $H$-coset.
\end{proof}

\begin{lemma}\label{NZ}  Let $S$ and $T$ be subsets of $\Z$  such that $|S|=3,$ $|T|=4$
and $|S+T|=7.$  Then $S$ is either a progression or  a
quasi-progression.
\end{lemma}

The proof is an easy exercise.

\begin{lemma}\label{lem:|T|=3}  Let $S$  be a $4$-separable
generating subset of an abelian  group $G$ of order $n$ such
that $0\in S$, $|S| = 3$ and $\kappa _4(S)=|S|=3.$ Assume moreover
that $\gcd (n,6)=1.$ Then $G$ is a cyclic group and $S$ is a
quasi-progression.\end{lemma}

\begin{proof*}
Put $S=\{0,x,y\}$. Let $0\in A$ be a $4$-atom of $S$. By Lemma
\ref{4at}, $|A|=4.$ Note that $A$ generates $G$ since otherwise
$|A+S|\ge 2|A|>|A|+|S|$.

We show first that every element of $ S\setminus \{ 0\}$ generates
$G$.  Suppose on the contrary that $x$ generates a proper subgroup
$K$ of $G$. Since $\gcd (|G|, 6)=1$ we have $\min\{ |H|, |G/K|\}\ge
5$.

Let $\phi$ denote the canonical morphism from $G$ onto $G/K$.
Decompose $A=A_1\cup\cdots \cup  A_j$, $j\ge 2$, modulo the subgroup
$K$ and assume that $0\in A_1$ and   $|A_1|\leq |A_i|$, $i\ge 2$.
Notice that
$$
|A+\{0,x\}|=\sum_{1\le i \le j} |A_i+\{0,x\}|\ge \sum_{1\le i \le j}
\min (|K|,|A_i|+1)\ge |A|+j.$$ On the other hand, since $\phi(S)$
generates $G/K$,  we have
$$|\phi (A+S)|=|\phi (A)+\phi (S)|\ge \min (|G/K|,|\phi (A)|+1)>|\phi(A)|=|\phi(A+\{ 0,x\})|.$$
Therefore,
$$
|A+S|\ge |A+\{0,x\}|+\min_i|A_i|\ge |A|+j+|A_1|,
$$
which implies $j=2$, $|A_1|=1$ and $|A_2|=3$. Now $A+S$ contains a
$K$--decomposition $ (A_2+\{0,x\}) \cup (A_2+y)$ involving only two
cosets. Thus $|A+S|\ge 1+|A_2|+ |A_2+\{0,x\}|\ge 8,$ a
contradiction. Hence each of $x$ and $y$ generate the cyclic group
$G=\zn$.

Since $|A|=4$ and $\gcd(|G|,6)=1$, we have $|A+\{ x,y\}|\ge |A|+1$. Assume
first that $|A+ \{ x,y\}|=|A|+1$. Then $A$ is an arithmetic
progression with difference $y-x$. But $0\in A$ and hence $y-x$ is
invertible since $A$ generates $G$. Without loss of generality we may assume
$A=\{0,1,2,3\}$. Now it comes easily that $S$ is a
quasi-progression, and the result holds.

 Suppose now that $|A+ \{ x,y\}|\ge |A|+2$. Since $|A+\{ 0,x,y\}|=|A|+3$, we may assume that
$|A\cap (A+x)|\geq 2$.

Now since $x$ is invertible in $G=\zn$, we may write, without loss
of generality, $S=\{0,1,t\}$ with $|A\cap (A+1)|\geq 2$. By
translating and multiplying by $-1$, we can also assume that $t\leq
(n+1)/2$ (notice that $\frac {n}{2}$ is not a unit if $n$ is even). Therefore $A$ can be represented by two pairs of
consecutive integers, and hence by a subset of $4$ integers
included in an interval of length $\leq (n+1)/2$. On the other
hand, one of the following two possibilities holds for $S$:

\begin{itemize}
\item  $S$ can be represented by a subset of an integral interval of length $\leq (n-3)/2$.
In that case the sum $A+S$ in $\zn$ has the same cardinality as
the sum $A+S$ in $\Z$, and we are done by Lemma \ref{NZ}.
\item   We have $t=(n-1)/2 $, in which case $S$ is included
in an arithmetic progression of length $4$ and difference $2^{-1}$
($2$ is invertible since $n$ is odd) and we are done. \  \rule{1mm}{2mm}
\end{itemize}
\end{proof*}

\section{Improving both the Theorems of Chowla and of Vosper}
\label{sec:chowlavosper}

Next we shall generalize Theorem \ref{th:kempermannis} to the case
when $|S+T|\leq |S|+|T|$. Our result is also a generalization to
abelian groups of Theorem \ref{th:HR}, i.e. the main result of
\cite{HR}. Let us state it
first under an isoperimetric formulation. Let us call a set
{\it quasi-periodic} if it can be obtained by deleting one element
from a periodic set.

\begin{theorem}\label{th:kempermannis+1}  Let $0\in S$ be a generating
  subset of a finite abelian
group $G$  with $\gcd(|G|,6)=1$ and $4\le |S|\le |G|-7$.
Assume $S$ to be $3$-separable and  $\kappa_3(S)= |S|.$

If every element of $S\setminus \{ 0\}$ has order at least
$|S|+1$, then either $S$ is a  quasi-progression or $S\setminus \{
0\}$ is quasi-periodic.
\end{theorem}

\begin{proof}
 Denote by $S^*=S\setminus \{ 0\}$. Let $0\in A$ be a $3$-atom of $S$.

\begin{claim}
The result holds if $A$ generates a
proper subgroup $K$ of $G$.
\end{claim}

Assume first that $A=K$. In this case $\kappa_3 (S)$ is a multiple
of $|A|$ and hence $|S|\geq |A|$. It follows  that $S\cap A=\{
0\}$, since otherwise $A$ would contain an element of order at
least $|S|+1$. Now $S+A$ is the disjoint union  $A\cup (S^*+A)$.
Hence $|S^*+A|=|S^*|+1$ so that $S^*$ is quasi-periodic and the
result holds.

Therefore we may assume that $A\neq K$. Decompose $S=S_1\cup\cdots
\cup  S_j$, $j\ge 2$, modulo the subgroup $K$. We may assume $0\in
S_1$ and $|S_1+A|\leq |S_2+A|\leq \cdots \leq  |S_j+A|.$ By Corollary
\ref{cor:2atom} and Lemma \ref{Sidon0},
\begin{equation}
  \label{eq:|S_i+A|}
  |S_i+A|\ge \min \{ |K|,|S_i|+|A|-1\}, \; i=1,\ldots j.
\end{equation}
It follows that $|S_i+A|= |K|$ for each $i\ge 2$. If $K+S=G$ then
$|S+A|\le |K+S|-3$ so that $|S_1+A|<|K|$. If $K+S\neq G$ then, by
the definition of $\kappa_3$, $|S+K|\ge |S|+|K|>|S|+|A|=|S+A|$ and
we also have also have $|S_1+A|<|K|$. It follows from
\eqref{eq:|S_i+A|} that
$|(S\setminus S_1)+A|\le |S\setminus S_1|+1$, which implies that
$S\setminus S_1$ is quasi-periodic.  In particular we have $|S|\ge
|K|$, which implies that $S_1=\{ 0\}$ and $S^*$ is quasi-periodic.
This completes the proof of the claim. 

We may therefore assume that $A$ generates $G$. We now consider
three cases.

{\it Case 1.\/}  $|A|=3$.

In that case $A$ is $4$-separable and $\kappa_4(A)\leq |A|$.
If $\kappa_4(A)<3$, then Theorem \ref{th:kempermannis} implies that
$A$ is a progression and thus $S$ is a quasi-progression. If
$\kappa_4(A)=3$, then
by Lemma \ref{lem:|T|=3},  $A$ is a
quasi-progression. By Lemma \ref{lem:transfer}, $S$ is a
quasi-progression.

{\it Case 2.\/}  $|A|=4$.

In that case $A$ is $3$-separable and $\kappa_3(A)\leq |A|$.
If $\kappa_3(A)<4$, then Theorem \ref{th:kempermannis}
implies that $A$ is a progression and thus $S$ is a quasi-progression.
If $\kappa_3(A)=4$ then
consider a $3$-atom $B$ of $A$ containing $0$. By
Lemma \ref{lem:|A|<=3}, $|B|=3.$ Observe that $B$ generates $G$,
otherwise $\gcd(|G|,6)=1$ implies $|A+B|\geq |A|+4>|A|+|B|.$
The set $B$ is $4$-separable and $\kappa_4(B)\leq |B|$.
If $\kappa_4(B)<3$ then $B$ is a progression
by Theorem \ref{th:kempermannis}, thus $A$, and hence also $S$ are
quasi-progressions. If $\kappa_4(B)=3$, then
by Lemma
\ref{lem:|T|=3},  $B$ is a quasi-progression. By Lemma
\ref{lem:transfer} applied twice we conclude that  $A$ and $S$ are
quasi-progressions.

{\it Case 3.\/} $|A|\ge 5$.

By Corollary \ref{cor:2atom},  for every
 $g\in G\setminus \{0\},$ we have
$$
|A\cap (A+g)|\leq 2.
$$

Let $B$ be a $2$-atom of $A$ containing $0$.
If $|B|=2$ then $|A|+|B|\ge |A+B|\ge
|A|+(|A|-2)$ which implies $|A|\le 4$. Hence we have $|B|\ge 3$.
By Theorem \ref{th:main}, $B$ is a subgroup of $G$. We have
$$
|A|+|B|\geq |A+B|\geq |A|+(|A|-2)+(|A|-4),
$$
and hence $|A|\geq \kappa_2(A)\geq |B|\geq  2|A|- 6.$ Since $\gcd
(|G|,6)=1$, these inequalities  force $|A|=|B|=5.$ It follows that
$|A+B|=2|B|$ and $A$ has a $B$-decomposition  $A=A_0\cup A_1.$ We
have $|A_i|\leq 3,$ since otherwise $|(A+x)\cap A|\geq 3$ for each
$x\in B\setminus \{ 0\}$, a contradiction. Without loss of
generality , we may assume $0\in A_0$, $|A_0|=3$ and $|A_1|=2.$

Decompose $S=S_0\cup \cdots \cup S_j$ modulo the subgroup $B$.

Assume first $S+B=G$. Let $V=\{i : S_i+B\not\subset S+A\}.$  By
Corollary \ref{chowla}, we have
$$
|A|=|A+S|-|S|\geq \sum _{i\in V} (|S_i+A_0|-|S_i|)\ge |V|(|A_0|-1).
$$
In particular $|V|\leq 2.$ Note that $|A+S|\le |G|-5$ since
otherwise $C=G\setminus (A+S)$ is a $3$-fragment of $-S$ and $-C$
a $3$-fragment of $S$ with $|C|<|A|$. Since
$|S_i+B|-|S_i+A_0|\le 2$, we have $|V|\ge 3$, a contradiction.

Assume now $S+B\neq G$.  Since $|A|>3$, Lemma \ref{lem:|A|<=3}
implies $|S|\ge 5=|B|$. Therefore $S\cap B=\{ 0\}$ since
every element of $S^*$ has order at least $|S|+1$.

Since $A$ generates $G$ we have $(A+S+B)\setminus (S+B)\neq
\emptyset$. Therefore, there is an $i$ such that $(A_1+S_i)\cap
(S+B)=\emptyset.$ Now $(S+A)\setminus S$ contains the disjoint
union $W= (A_0\setminus \{ 0\}) \cup  (S_i+A_1).$ But $|W|\geq
|A|-1.$ It follows easily that $S\setminus \{ 0\}$ is periodic  or
quasi-periodic. The first possibility is excluded by the condition
$\kappa _3(S)=|S|.$
\end{proof}

Theorem \ref{th:kempermannis+1} translates into a
characterization of subsets $S$ and $T$ such that $|S+T|\le
|S|+|T|$ under some Chowla-type conditions. This was our final
goal in this paper.

\begin{theorem}\label{th:kemperman+1}
Let $G$ be  a finite abelian group  with $\gcd(|G|,6)=1$.

Let $0\in S$ be a generating subset of $G$ such that
$|S|\ge 4$ and every element in $S^*$ has order at least $|S|+1$.
Let $Q$ be a maximal
subgroup   such that $|S^*+Q|-|S^*|\leq 1$  and let
$\sigma:G\rightarrow G/Q$ denotes the canonical projection.

Let $T$ be a subset of $G$ such that $|T|\ge 3$ and suppose that
$|S+T|= |S|+|T|\le |G|-4$. Then the following holds:
\begin{itemize}
\item If $Q=\{0\}$ then  $S$ and $T$ are progressions or
quasi-progressions with the same difference.
\item If $Q\neq \{0\}$,
$|\sigma(T)|\geq 2$, and $|\sigma(S+T)|<|G|/|Q|-1$ then
$\sigma(S)$ and $\sigma(T)$ are arithmetic progressions with the same
 difference.
Moreover we have
$|(T\setminus T_1)+Q|\leq |T\setminus T_1|+1$ where $T_1$
is a subset of $T$ such that $\sigma(T_1)$ is a single, extremal
element of the progression $\sigma(T)$.
Furthermore  if $0$ is not an extremal
element of the progression $\sigma(S)$, then $|T+Q|\leq |T|+1$.
\end{itemize}
\end{theorem}

\begin{proof}

\noindent
{\bf Case 1: $Q=\{0\}$.}\\
The conditions $|S|+|T|=|S+T|\le |G|-4$ and $|T|\ge 3$ imply that
$S$ is $3$-separable and that $\kappa _3(S)\leq |S|.$ By Theorems
\ref{th:kempermannis} and \ref{th:kempermannis+1}, $S$ is an
arithmetic progression or quasi-progression. By Lemma
\ref{lem:transfer} it follows that $T$ is an arithmetic
progression or  quasi-progression with the same difference.

\medskip

\noindent
{\bf Case 2: $Q\neq\{0\}$.}\\
Since $|S|\ge q=|Q|$ and each element in $S$ has order at least
$|S|+1$, we have $S\cap Q=\{ 0\}$. In particular, $\sigma
(S)^*=\sigma (S^*)$.
Note that each element in $\sigma (S)^*$ has
order at least $(|S|+1)/q\ge |\sigma (S)|-1+1/q$.

\noindent
Corollary \ref{chowla} implies
\begin{equation}\label{eq:withcor8}
  |\sigma(S)+\sigma(T)|\geq |\sigma (S)|+|\sigma (T)|-1.
\end{equation}
First notice that
\begin{equation}\label{eq:Qperiodic}
  \Sigma = (S^*+T)\setminus (T+Q) \;\;\text{is $Q$-periodic.}
\end{equation}
This holds clearly if  $S^*$ is $Q$-periodic.
So we may assume   $|S^*+Q|-|S^*|=1$.
Let us then denote by $S_1$
the unique subset of $S$ of size $|Q|-1$ in the decomposition of
$S$ modulo $Q$. If $\Sigma$ is not $Q$-periodic then
some $Q$-coset must have a trace $U$ of size $|Q|-1$ on the set
$\Sigma$, and we have
$U=S_1+T'$ where  $T' = (a+Q)\cap T$, for some $a$.
Since
$|S_1|=|Q|-1$ we must have $|T'|=1$. Note also
that $\sigma(S_1)+\sigma(T')$ cannot be obtained in any other way
as a sum of an element of $\sigma(S)$ and of an element of $\sigma(T)$,
therefore $(S_1+T') \cap (S+(T\setminus T')) = \emptyset$, hence
$|S+(T\setminus T')| < |S|+ |T\setminus T'| -1$, but this contradicts
$\kappa_1(S)=|S|-1$ (Corollary \ref{chowla}) and proves \eqref{eq:Qperiodic}.

We therefore have:
\begin{equation}
  \label{eq:noname}
|\sigma(S)+\sigma(T)| = |\sigma(S)|+|\sigma(T)|-1.
\end{equation}
otherwise \eqref{eq:withcor8} implies
$|S+T|\geq |T| + |\sigma(S)||Q| >|T|+|S|$ a contradiction.

By our assumptions, $Q$ is a maximal
subgroup   such that $|S^*+Q|-|S^*|\leq 1$. This is easily seen to imply
that $\sigma (S^*)$ is not periodic.
 Moreover, each element in
$\sigma (S)^*$ has order at least $(|S|+1)/q\ge |\sigma
(S)|-1+1/q$. Then,
by Proposition \ref{fragments}, $\sigma (S)$ and $\sigma (T)$ are
arithmetic progressions with the same common difference $d$.
Since $-d$ is also a difference of $\sigma (S)$ and $\sigma (T)$,
we may assume without loss
of generality that  the terminal element $u$ of $\sigma(S)$ is not
 $0$. Therefore if we set $S'=\sigma^{-1}(u)\cap S$ we have
$|S'|\geq |Q|-1$. Let us suppose, without loss of generality, that
the initial element of $\sigma(T)$ is $0$.

To conclude we prove the following:

\begin{claim}
  Let  $B\subset G$ be  such that $\sigma(B)$ is an arithmetic
  progression of difference $d$ and initial element $0$,
  and with $|\sigma(B)|\leq |\sigma(T)|$.
  Let $B_1=\sigma^{-1}(0)\cap B$ and set
    $B_2=B\setminus B_1$. Suppose  that $|S+B|=|S|+|B|-\varepsilon$ where
  $\varepsilon$ equals $0$ or $1$. Then $|B_2+Q|\leq |B_2|+1-\varepsilon$.
\end{claim}

The claim is proved by induction on $t=|\sigma(B)|$.
If $|\sigma(B)|=1$ then $B_2=\emptyset$ and there is nothing to prove.
If $|\sigma(B)|\geq 2$, then let $b$ be the terminal element of
$\sigma(B)$ and set $B'=\sigma^{-1}(b)\cap B$. Then, since
$\kappa_1(S)=|S|-1$:
\begin{equation}
|S|+|B|-|B'|-1\leq |S+(B\setminus B')|
\le |S+B|-|S'+B'|=|S|+|B|-\varepsilon -|S'+B'|.
\label{eqFINAL}
\end{equation}
{F}rom \eqref{eqFINAL} we obtain, since $|S'|\geq |Q|-1$,
 $|B'|\ge |S'+B'|-1\ge |Q|-2\ge 5-2=3.$
Hence, since $|S'|+|B'|>|Q|$, we have $S'+B'=Q$.
Applying \eqref{eqFINAL} again we also have $|S'+B'|\leq
|B'|+1-\varepsilon$,
hence
\begin{equation}
  \label{eq:new}
  |B'|\geq |Q|-1+\varepsilon
\end{equation}
By \eqref{eqFINAL} we have
$$|S+(B\setminus B')|\leq |S|+|B\setminus B'|+(|B'|-|S'+B'|)
-\varepsilon .$$
Now either $|B'|-|S'+B'|=0$
and the result holds by the induction hypothesis applied to
$B\setminus B'$: or $|B'|-|S'+B'|=-1$. But in this case
\eqref{eq:new} implies $\varepsilon=0$, and the result again holds
by applying the induction hypothesis to $B\setminus B'$ with
$|S+(B\setminus B')|\leq |S|+|B\setminus B'|-1$.
This proves the claim and the theorem with $T_1=\sigma^{-1}(0)\cap T$.

To complete the proof, consider the case when $0$ is not an
extremity of $\sigma(S)$. Then the claim applied to $S$ and $T$ by
permuting initial and terminal elements of $\sigma(S)$ and
$\sigma(T)$ gives $|T_1|\geq |Q|-1$, so that every non-empty
intersection of $T$ with a coset of $Q$ has cardinality at least
$|Q|-1$. In particular, $S^*+T$ is $Q$--periodic. Since $0$ is not
an extremal element of $\sigma(S)$ and $|\sigma(T)|\geq 2$ we obtain
$\sigma(S^*)+\sigma (T)=\sigma (S)+\sigma (T)$ so that $S+T$ is also
$Q$-periodic. By \eqref{eq:noname} we have $|T+Q|\leq |T|+1$.
\end{proof}

\noindent
{\bf Remark.} One may wonder what happens if we remove from
Theorem \ref{th:kemperman+1} the hypothesis
$|\sigma(S+T)|<|G|/|Q|-1$. Then the sets $\sigma(S)$ and $\sigma(T)$ are
not necessarily arithmetic progressions any more. However, one may show
that there again exists $T_1\subset T$, such that $|\sigma(T_1)|\leq 1$ and
$T\setminus T_1$ is $Q$-periodic or $Q$-quasi-periodic.
We leave out the details.

\section*{Acknowledgements} The authors are grateful to an anonymous
referee for his/her careful reading and very valuable suggestions
and remarks.

\end{document}